\documentclass{ifacconf}

\usepackage{graphicx}      				
\usepackage{natbib}        				
\usepackage{amsmath,amssymb,amsfonts}	
\usepackage{siunitx} 
\usepackage{color}

\newcommand{\Int}{\mathsf{int}}

\newcommand{\U}{\mathcal{U}}
\newcommand{\M}{\mathcal{M}}
\newcommand{\DMM}{\left[\partial \mathcal{M}\right]_-}

\newcommand{\A}{\mathcal{A}}
\newcommand{\Acomp}{\mathcal{A}^{\textsf{C}}}
\newcommand{\DAM}{\left[\partial \mathcal{A}\right]_-}
\newcommand{\DAO}{\left[\partial \mathcal{A}\right]_0}

\newtheorem{definition}{Definition}

\begin{document}
\begin{frontmatter}

\title{Maintaining Hard Infection Caps in Epidemics via the Theory of Barriers\thanksref{footnoteinfo}} 

\thanks[footnoteinfo]{This project has received funding from the European Social Fund (ESF).}

\author[First]{Willem Esterhuizen}
\author[First]{Tim Aschenbruck}
\author[Second]{Jean L\'{e}vine}
\author[First]{Stefan Streif}

\address[First]{Automatic Control and System Dynamics Laboratory, Technische Universit\"{a}t Chemnitz, 09111 Chemnitz, Germany (e-mail: tim.aschenbruck@etit.tu-chemnitz.de; willem.esterhuizen@etit.tu-chemnitz.de; stefan.streif@etit.tu-chemnitz.de).}
\address[Second]{CAS, Math\'{e}matiques et Syst\`{e}mes, Mines-ParisTech, PSL, 60, Boulevard Saint-Michel, 75006 Paris, France (e-mail: jean.levine@mines-paristech.fr).}

\begin{abstract}
	Research in epidemiology often focusses on designing interventions that result in the number of infected individuals asymptotically approaching zero, without considering that this number may peak at high values during transients. Recent research has shown that a set-based approach could be used to address the problem, and we build on this idea by applying the theory of barriers to construct admissible and invariant sets for an epidemic model. We describe how these sets may be used to choose intervention strategies that maintain infection caps during epidemics. We also derive algebraic conditions of the model parameters that classify a system as being either comfortable, comfortable-viable, viable, or desperate.
\end{abstract}

\begin{keyword}
State Constraints, Epidemic, Infection Cap, Admissible Set, Maximal Robust Positively Invariant Set (MRPI), Barrier
\end{keyword}

\end{frontmatter}

\section{Introduction}

Epidemics have classically been modelled by so-called \emph{compartmental models}, where the state variables represent the proportion of individuals that belong to a certain class (susceptible, infected, recovered, exposed, vaccinated, etc.) with their interaction governed by ordinary differential equations (ODEs). See for example \cite{brauer_2001_mathematical_models_population} and \cite{Hethcote_maths_of_infectuous_diseases_2000}. Research in the management of epidemics have usually focused on the stability properties of these ODEs under (stationary) interventions.

Control theory has been applied in this classical ODE setting, with researchers often analysing optimal control problems that optimise a trade-off between peak proportion of infected members and the cost of intervention, see for example \cite{Culshaw2004}, \cite{Hansen2011}, \cite{HETHCOTE1973365}, \cite{Kirschner1997}. There has also been much recent interest in the modelling of epidemics as complex networks, see \cite{CASTELLANO_2015}, and \cite{Nowzari_2016} for a survey of control-theoretic approaches to manage epidemics.

The mentioned papers that analyse optimal control problems aim at minimising the infected population by imposing a relevant cost (\emph{soft} constraint). Some work has been done on analysing the size of disease outbreaks, see \cite{miller_2008} for a stochastic-model setting. However, to our knowledge, the only papers that address maintaining \emph{hard} constraints on the infected population are those by \cite{deLara2016viableMOSQUITOS,Delara_viab_mosquitos_2019}, which we briefly discuss.

\cite{deLara2016viableMOSQUITOS} first proposed the use of \emph{viability theory}, \cite{aubin2009viability}, to address the problem of maintaining a hard infection cap. They found the \emph{viability kernel}  for the Ross-Macdonald model, \cite{anderson1991BOOKRossMacDonaldModel}, which describes the dynamic behaviour of a vector-borne disease (a vector is an agent, such as a mosquito, that transmits disease-causing pathogens). They defined three cases that may occur, namely the \emph{comfortable}, \emph{viable} and \emph{desperate} cases, which correspond to inequalities on the system parameters. In particular, the viable case corresponds to when there exists a nontrivial\footnote{neither empty nor equal to the constraint set}viability kernel. In this case, the authors describe a special part of the boundary of the set, arguing that it is given by a solution of the system obtained with a particular input.

In this paper we use the theory of \emph{barriers}, as developed in \cite{DeDona_Levine2013barriers} and \cite{Ester_Asch_Streif_2019}, to analyse the Ross-Macdonald model subjected to state and input constraints. We derive necessary conditions that must be satisfied by the \emph{admissible set}\footnote{The admissible set is the set of initial conditions such that there exists at least one admissible control for which the corresponding integral curve satisfies the state constraints for all times (see Definition~\ref{adm:def}). It is called the viability kernel in viability theory} or the \emph{maximal robust positively invariant set} (MRPI)\footnote{The MRPI is the set of initial conditions for which the state constraints are always satisfied, regardless of the admissible input (see Definition~\ref{MRPI:def}).} for the system. When these sets are non trivial it is easy to find them: one merely needs to solve a differential equation over some interval of time. Intervention strategies that maintain infection caps can then be determined depending on the set that the state belongs to. We elaborate on this point at the end of Section~\ref{sec:2}.

Our contributions are as follows:
\begin{itemize}
	\item We analyse the system with bounds on the proportion of humans \emph{and} mosquitoes. If the mentioned sets are non trivial, then they also give information on whether limits on the proportion of infected mosquitoes can be maintained.
	\item In addition to defining the comfortable, viable and desperate cases for epidemic models, we introduce the \emph{comfortable-viable case}, which occurs when the system has a nontrivial MRPI.
	\item Using the main result from the theory of barriers, we are able to state inequalities on the system parameters under which each case occurs.
	\item We reproduce the conditions derived by \cite{deLara2016viableMOSQUITOS} for the comfortable, viable and desperate cases, showing that they are obtained by setting the mosquito cap, $\bar{x}_1$, equal to one.
	\item It is interesting to note that similar to \cite{deLara2016viableMOSQUITOS}, we describe certain parts of the boundaries of the sets (the ``barriers'') as being made up of special system curves, not only for the admissible set but also for the MRPI: they satisfy a maximum/minimum-like principle.
	\item The comparisons between the admissible set and the MRPI may lead to interesting conclusions concerning the epidemic management, thus building on the results of \cite{deLara2016viableMOSQUITOS}.
\end{itemize}

The outline of the paper is as follows. In Section~\ref{sec:2} we present a summary of the theory of barriers along with other mathematical preliminaries, such as our notation. In Sections~\ref{sec_RM} we present an in-depth analysis of the Ross-Macdonald model, the main contributions being the conditions that determine the different system cases. In Section~\ref{sec_examples} we present a numerical example. We discuss how the obtained sets may be used in the management of epidemics in Section~\ref{sec_discussion}, and conclude the paper in Section~\ref{sec_conclusion}.

\section{Summary of the Theory of Barriers}\label{sec:2}
We briefly summarise the results as presented in \cite{DeDona_Levine2013barriers} and \cite{Ester_Asch_Streif_2019}.

We consider the following constrained nonlinear system:
\begin{align}
\dot{x}(t) & = f(x(t),u(t)),\label{sys_eq_1}\\
x(t_0) & = x_0,\label{sys_eq_2}\\
u & \in \mathcal{U},\label{sys_eq_3}\\
g_i(x(t)) &\leq 0,\forall t\in[t_0,\infty[,\,\,i=1,2,\dots,p,\label{sys_eq_4}
\end{align}
where $x(t)\in\mathbb{R}^n$ denotes the state, $u(t)\in\mathbb{R}^m$ denotes the input, and $\mathcal{U}$ is the set of Lebesgue measurable functions that map the interval $[t_0,\infty[$ into a compact and convex set $U\subset \mathbb{R}^m$. The initial condition is specified by $x_0$ and the $g_i$'s are the constraint functions. Briefly, the assumptions are that the functions $f$ and $g_i$ for $i=1,\dots,p$ are $C^2$ with respect to their arguments on appropriate open sets; that all solutions of the system remain bounded on finite intervals; and that the set $\{f(x,u):u\in U\}$ is convex for all $x$. See  \cite{DeDona_Levine2013barriers} and \cite{Ester_Asch_Streif_2019} for more details. It is easy to verify that the system under study in this paper satisfy these assumptions.

We will use the notation $x^{\bar{u}}$ to denote a solution to \eqref{sys_eq_1} with input $\bar{u}\in\mathcal{U}$, where the initial condition is clear from context. We denote $x^{(u,x_0,t_0)}(t)$ as the solution of \eqref{sys_eq_1} at time $t$, with initial condition \eqref{sys_eq_2} at $t_0$ and input \eqref{sys_eq_3}. We introduce the following sets:
\begin{align}
G\triangleq& \{x : g_i(x)\leq 0, \quad  \forall i\in\{1,2,...,p\}\},\label{def_G}\\
G_-\triangleq& \{x : g_i(x) <0, \quad  \forall i\in\{1,2,...,p\}\},\nonumber\\
G_0\triangleq& \{x \in G: \exists i\in\{1,2,...,p\} \;\; \text{s.t.} \;\; g_i(x) =0\},\nonumber
\end{align}
and by $\mathbb{I}(x)=\{i: g_i(x) = 0\}$, the indices of active constraints at $x$.
Furthermore, we denote by $L_fg(x,u)\triangleq Dg(x)f(x,u)$ the Lie derivative of a continuously differentiable function $g:\mathbb{R}^n \rightarrow \mathbb{R}$ with respect to $f(.,u)$ at the point $x$.
Given a set $S$, its boundary is denoted by $\partial S$ and $S^{\textsf{C}}$ denotes its complement. By $\mathbb{R}_{\geq 0}$ we refer to the set of nonnegative real numbers. 

\begin{definition}\label{adm:def}
	The \emph{admissible set} of the system \eqref{sys_eq_1}-\eqref{sys_eq_4}, denoted by $\A$, is the set of initial states for which there exists a $u\in \U$ such that the corresponding solution to \eqref{sys_eq_1}-\eqref{sys_eq_2} satisfies the constraints \eqref{sys_eq_4} for all future time.
	\begin{align*}
	\mathcal{A} \triangleq \left\{ x_0 \in G: \exists u\in\mathcal{U}, \;\;  x^{(u,x_0,t_0)}(t)\in G \;\; \forall t \in [t_0,\infty[ \;\right\}.
	\end{align*}
\end{definition}
\begin{definition}\label{RPI:def}
	A set $\Omega\subset G$ is said to be a \emph{robust positively invariant set} (RPI) of the system \eqref{sys_eq_1}-\eqref{sys_eq_3} provided that $x^{(u,x_0,t_0)}(t)\in\Omega$ for all $t\in[t_0,\infty[$, for all $x_0\in\Omega$ and for all $u\in\U$.
\end{definition}
\begin{definition}\label{MRPI:def}
	The \emph{maximal robust positively invariant set} (MRPI) of the system \eqref{sys_eq_1}-\eqref{sys_eq_4} contained in $G$, is the union of all RPIs that are subsets of $G$. Equivalently\footnote{shown in \cite{Ester_Asch_Streif_2019}, Proposition~2}
	\begin{align*}
	\mathcal{M} \triangleq \left\{ x_0\in G: x^{(u,x_0,t_0)}(t)\in G, \,\,\forall u\in\mathcal{U}, \;\;  \forall t \in [t_0,\infty[ \;\right\}.
	\end{align*}
\end{definition}

We introduce the sets:
\[
	\DAM \triangleq  \partial\mathcal{A}\cap G_-, \quad	\DMM  \triangleq \partial\mathcal{M}\cap G_-,
\]
where $\DAM$ and $\DMM$ are called the \emph{barrier} and \emph{invariance barrier}, respectively. These parts of the boundaries of the sets have the property that for every point $\bar{x}\in\DAM$ (resp. $\bar{x}\in\DMM$) there exists an input, $\bar{u}\in \mathcal{U}$, such that the resulting integral curve runs along $\DAM$, (resp. $\DMM$), and satisfies a minimum/maximum-like principle. Moreover, if the curves eventually intersect the boundary of the constraint set, then this must happen in a tangential manner. These facts are summarised in the following theorem.
\begin{thm}\label{thm1}
	Under the assumptions in \cite{DeDona_Levine2013barriers} and \cite{Ester_Asch_Streif_2019}, every integral curve $x^{\bar{u}}$ on $\DAM$ (resp. $\DMM$) and the corresponding input function $\bar{u}$ satisfy the following necessary conditions. There exists a nonzero absolutely continuous maximal solution $\lambda^{\bar{u}}$ to the adjoint equation:
	\begin{align}
	&\dot{\lambda}^{\bar{u}}(t) = -\left( \frac{\partial f}{\partial x}(x^{\bar{u}}(t),\bar{u}(t)) \right)^T \lambda^{\bar{u}}(t),\nonumber
	\end{align}
	such that
	\begin{align}		
	\min_{u\in U}&\{\lambda^{\bar{u}}(t)^Tf(x^{\bar{u}}(t),u) \} \nonumber\\ 
	&= \lambda^{\bar{u}}(t)^T f(x^{\bar{u}}(t),\bar{u}(t)) = 0 \label{eq_Hamil_A} \\
	\Big(\Big.\text{resp.}\quad
	\max_{u\in U}&\{\lambda^{\bar{u}}(t)^Tf(x^{\bar{u}}(t),u) \} \nonumber\\ &= \lambda^{\bar{u}}(t)^T f(x^{\bar{u}}(t),\bar{u}(t)) = 0\label{eq_Ham_max_M}
	\Big.\Big),
	\end{align}
	for almost all $t$. Moreover, if $x^{\bar{u}}$ intersects $G_0$ in finite time, we have:
	\begin{equation}
	\lambda^{\bar{u}}(\bar{t}) = (Dg_{i^*}(z))^T,\label{eq_adj}
	\end{equation}
	where
	\begin{align}
	&\min_{u\in U}\max_{i\in\mathbb{I}(z)} L_fg_i(z,u) = L_fg_{i^*}(z,\bar{u}(\bar{t})) = 0 \label{thm1_ult_tan_A} \\
	\Big(\Big.\text{resp.}\quad
	&\max_{u\in U}\max_{i\in\mathbb{I}(z)} L_fg_i(z,u) = L_fg_{i^*}(z,\bar{u}(\bar{t})) = 0 \label{thm1_ult_tan_M}
	\Big.\Big),
	\end{align}
	$\bar{t}$ denotes the time at which $G_0$ is reached, and $z \triangleq x^{(\bar{u},\bar{x},t_0)}(\bar{t})\in G_0$.
\end{thm}
These necessary conditions may be used to construct the sets $\A$ and $\M$ as follows. First, identify so-called points of \emph{ultimate tangentiality} (located on $G_0$) via the conditions \eqref{thm1_ult_tan_A} or \eqref{thm1_ult_tan_M}. To be less verbose, in the sequel we will also refer to these points as \emph{tangent points}. Then, consider the system's Hamiltonian, $H(x,\lambda,u) \triangleq \lambda^Tf(x,u)$, and determine the input realization $\bar{u}$ corresponding to the integral curve along $\DAM$ or $\DMM$ using the \emph{Hamiltonian minimisation/maximisation condition}, \eqref{eq_Hamil_A} or \eqref{eq_Ham_max_M}. Finally, find the integral curves along $\DAM$ or $\DMM$ by integrating the system and adjoint dynamics backwards from the tangent points.

Having found the sets $\A$ and $\M$ for a constrained system, information can be determined from the location of the system's current state. If the state is located in $\M$, then the system is \emph{safe}: the state constraints are guaranteed to be maintained for all time, regardless of the input function. If the state is located in $\A$, then the system is \emph{potentially safe}: it is possible to specify an input function such that the state constraints are always satisfied. If the state is located in $\Acomp$, the system is \emph{unsafe}: the state constraints are guaranteed to be violated sometime later, regardless of the input function chosen. We will say more on how these facts translate into epidemic management strategies in Section~\ref{sec_discussion}.

We note that the conditions of the theorem are necessary, and so it may be that certain parts of the obtained curves need to be ignored. Thus, we will refer to integral curves obtained via the necessary conditions as \emph{candidate barrier} and \emph{candidate invariance barrier} curves.

We also note that the adjoint $\lambda^{\bar{u}}$ associated with a barrier/invariance barrier curve $x^{\bar{u}}$ is the normal of a hyperplane that evolves along the curve. From \eqref{eq_Hamil_A} and \eqref{eq_Ham_max_M} we see that $\lambda^{\bar{u}}$ is perpendicular to the curve, a fact we will use in Lemma~\ref{lemma_2}.

\section{Analysis of Ross-Macdonald model}\label{sec_RM}

In this section we carry out an analysis of a constrained Ross-Macdonald model, using the conditions of Theorem~\ref{thm1}.

\subsection{Ross-Macdonald model}
We consider the following constrained model \cite[p. 394]{anderson1991BOOKRossMacDonaldModel} for a mosquito-borne disease:
\begin{align}\label{eq_RossMacdonald}
\begin{split}
\dot{x}_1 &= A_m x_2(1-x_1) - u x_1, \\
\dot{x}_2 &= A_h x_1 (1 - x_2) - \gamma x_2,\\
x_1(t) &\in [0,\bar{x}_1],\\
x_2(t) &\in [0,\bar{x}_2],\\
u(t) &\in U \triangleq \left[u_{\min}, u_{\max}\right],
\end{split}
\end{align}
where $x_1$ is the proportion of infected mosquitoes, $x_2$ is the proportion of infected humans, and $u$ is the control which encapsulates the natural mortality rate of the mosquitoes, as well as the effects of fumigation. We let $0 < u_{\min} < u_{\max}$. The constants $A_m \triangleq a p_m \geq 0$ and $A_h \triangleq a p_h \frac{N_m}{N_h} \geq 0$, where $a$ denotes the biting rate, $p_m$ is the probability of a mosquito being infected after
biting an infected human, $p_h$ denotes the probability of infection of a susceptible human
after being bitten by infected mosquito, $\frac{N_m}{N_h}$ is the ratio of female mosquitoes to humans and $\gamma$ denotes the human recovery rate. The constants satisfy: $a \geq 0$, $\frac{N_m}{N_h}~\geq~0$, $0 \leq p_h \leq 1$, $0 \leq p_m \leq 1$, $\gamma \geq 0$.
The bounds $0 < \overline{x}_1 \leq 1$ $0 < \overline{x}_2 \leq 1$ denote the maximum tolerated proportion of infected mosquitoes and humans, respectively.

\subsection{Points of ultimate tangentiality}\label{subsec_ult_tan}
We label the state constraints as follows: 
$g_1(x) = x_1 - \bar{x}_1$, $g_2(x) = - x_1$, $g_3(x) = x_2 - \overline{x}_2$, $g_4(x) = - x_2$, and the $i$-th tangent point on the constraint $g_i$ as $z^i \triangleq (z^i_1,z^i_2)^T$.
Concentrating on $g_2$ and $g_4$, and invoking \eqref{thm1_ult_tan_M}, we get:
\begin{align*}
\max_{u\in U } \left\{L_fg_2(z,u) \right\} &= -A_mz^2_2 < 0,\\
\max_{u\in U } \left\{L_fg_4(z,u) \right\} &= -A_hz^4_1 < 0.
\end{align*}
We see that at all points for which $g_2$ and $g_4$ are active the state evolves (forwards in time) into the interior of the constrained state-space regardless of the chosen input. The constraint functions $g_1$ and $g_3$ are more interesting. We see that:
\[
\max_{u\in U} \left\{L_fg_1(z,u)\right\} = A_mz^1_2(1 - \bar{x}_1) - u_{\min}\bar{x}_1,
\]
and setting this to zero, we get a tangent point associated with the MRPI at $z^1 = (\bar{x}_1, \frac{u_{\min}\bar{x}_1}{A_m(1 - \bar{x}_1)})^T$. Note that $\frac{u_{\min}\bar{x}_1}{A_m(1 - \bar{x}_1)}\in[0,\infty[$ for all $\bar{x}_1\in[0,1]$. Thus, $z^1\in \{\bar{x}_1\}\times [0,\bar{x}_2[$ if and only if:
\begin{equation}
	\bar{x}_1 < \frac{A_m\bar{x}_2}{A_m\bar{x}_2 + u_{\min}}.\label{existence_ult_tan_M_g1}
\end{equation}
Similarly, $\min_{u\in U} L_f g_1(z,u) = 0$ gives us 
$$z^1 = (\bar{x}_1, \frac{u_{\max}\bar{x}_1}{A_m(1 - \bar{x}_1)})^T$$ 
associated with the admissible set, which is located on $\{\bar{x}_1\}\times [0,\bar{x}_2[$ if and only if:
\begin{equation}
\bar{x}_1 < \frac{A_m\bar{x}_2}{A_m\bar{x}_2 + u_{\max}}.\label{existence_ult_tan_A_g1}
\end{equation}
Concentrating on $g_3$, we get:
\begin{align}
	L_fg_3(z,u) = A_hz^3_1(1 - \overline{x}_2) - \gamma \overline{x}_2,\label{eq_LFG3}
\end{align}
which is independent of the input, giving the tangent point $z^3 = (\frac{\gamma\overline{x}_2}{A_h (1 - \overline{x}_2)}, \overline{x}_2)^T$ for both the admissible set and the MRPI. We note that $\frac{\gamma\overline{x}_2}{A_h (1 - \overline{x}_2)}\in[0,\infty[$ for all $\bar{x}_2\in[0,1]$. Thus, $z^3 \in [0,\bar{x}_1[\times \{\overline{x}_2\}$ if and only if:
\begin{equation}
	\bar{x}_2 < \frac{A_h\bar{x}_1}{A_h\bar{x}_1 + \gamma}.\label{existence_ult_tan}
\end{equation}

\subsection{Barrier curves}\label{subsec_input_realisation}
Invoking \eqref{eq_Hamil_A} and \eqref{eq_Ham_max_M}, we obtain:
\begin{align*}
\bar{u}(t) = 
\begin{split}
\begin{cases}
u_{\max} & \text{if} \quad \lambda_1(t) \geq 0\\
u_{\min} & \text{if} \quad \lambda_1(t) < 0,
\end{cases}
\end{split}
\end{align*}
for the input associated with the admissible set, and
\begin{align*}
\bar{u}(t) = 
\begin{split}
\begin{cases}
u_{\min} & \text{if} \quad \lambda_1(t) \geq 0\\
u_{\max} & \text{if} \quad \lambda_1(t) < 0,
\end{cases}
\end{split}
\end{align*}
for the input associated with the MRPI. The adjoint equation is:
\begin{align} \label{eq_adj_ross_macdonald}
\dot{\lambda}=&
\begin{pmatrix}
A_m x_2+ u & -A_h(1-x_2)\\
-A_m(1-x_1) & A_h x_1+\gamma
\end{pmatrix}
\lambda,
\end{align}
with $\lambda(\bar{t}) = (1,0)^T$ associated with $z_1$, and $\lambda(\bar{t}) = (0,1)^T$ associated with $z_3$.

\begin{figure}[h]
	\begin{center}
		\includegraphics[width=\linewidth,height=6cm,keepaspectratio]{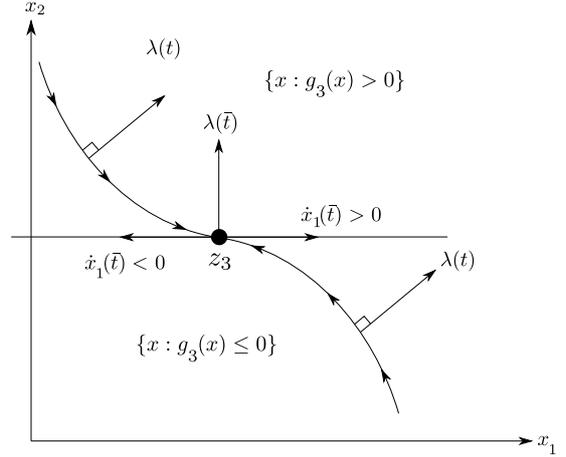} 
		\caption{Clarification of the proof of Lemma~\ref{lemma_2}: knowing that $\lambda_1(t)>0$ before the barrier intersects $z_3$, we must have $\dot{x}_1(\bar{t}) <0$.} 
		\label{fig:lem_2}                             
	\end{center}                        
\end{figure}

We now state some lemmas that summarise important aspects of the system's barrier curves. We will use these results in Proposition~\ref{prop_cases}.
\begin{lem}\label{lemma_2}
	Consider the  system \eqref{eq_RossMacdonald} with $\bar{x}_1 < 1$ and $\bar{x}_2 < 1$. There exists a candidate barrier curve associated with $\A$ (resp. $\M$), partly contained in $G_-$ and ending at:
	\begin{itemize}
		\item $z^1 = (\bar{x}_1, \frac{u_{\max}\bar{x}_1}{A_m(1 - \bar{x}_1)})^T$ $\left(\mathrm{resp.}\,\,z^1 = \bar{x}_1, \frac{u_{\min}\bar{x}_1}{A_m(1 - \bar{x}_1)})^T\right)$ if and only if
	\begin{align}
		A_h(A_m + u_{\max})\bar{x}_1 + \gamma u_{\max} & > A_mA_h,\label{viable_case_z1}\\
		(\textrm{resp.}\quad A_h(A_m + u_{\min})\bar{x}_1 + \gamma u_{\min} &> A_mA_h \label{MRPI_case_z1}),
	\end{align}
	\item $z^3 = (\frac{\gamma\overline{x}_2}{A_h(1-\overline{x}_2)}, \overline{x}_2)$ if and only if
	\begin{align}
		A_m(A_h + \gamma)\bar{x}_2 + \gamma u_{\max} &> A_mA_h,\label{viable_case}\\
		(\textrm{resp.}\quad 	A_m(A_h + \gamma)\bar{x}_2 + \gamma u_{\min} &> A_mA_h \label{unconquerable_case_MRPI}).
	\end{align}
	\end{itemize}
\end{lem} 
\begin{pf}
	At the tangent point, $z^1$, we have $\lambda(\bar{t}) = (1,0)^T$, and thus $\bar{u}(t) = u_{\max}$ on an interval before $\bar{t}$ for the input associated with $\A$, and $\bar{u}(t) = u_{\min}$ for the input associated with $\M$.	From \eqref{eq_adj_ross_macdonald} we see that $\dot{\lambda}_2(\bar{t}) < 0$, implying that $\lambda_2(t) > 0$ on an interval before $\bar{t}$. Therefore, the integral curve will evolve backwards into $G_-$ if and only if $\dot{x}_2(\bar{t}) < 0$. Substituting the value for $z^1$ we arrive at the statements \eqref{viable_case_z1} and \eqref{MRPI_case_z1}. The cases for $z^3$ follow similar arguments: we have $\lambda(\bar{t}) = (0,1)^T$, and we see that $\dot{\lambda}_1(\bar{t}) < 0$, implying that $\lambda_1(t) \geq 0$ on an interval before $\bar{t}$. Thus, $\bar{u}(t) = u_{\max}$ on an interval before $\bar{t}$ for the input associated with $\A$, and $\bar{u}(t) = u_{\min}$ for the input associated with $\M$. Because $\lambda_1(t) > 0$ before $\bar{t}$, the integral curve will evolve backwards into $G_-$ if and only if $\dot{x}_1(\bar{t}) < 0$, see Figure~\ref{fig:lem_2} for clarification. Substituting the value for $z^3$ we arrive at the statements. This completes the proof.
\end{pf}

\begin{lem}\label{lemma_3}
	Suppose the system \eqref{eq_RossMacdonald} has candidate barrier curves associated with $\A$ (resp. $\M$) partly contained in $G_-$ and ending at a tangent point, as in Lemma~\ref{lemma_2}. Then, the inequalities \eqref{existence_ult_tan_A_g1} (resp. \eqref{existence_ult_tan_M_g1}) and \eqref{existence_ult_tan} cannot both hold.
\end{lem}
\begin{pf}
	Suppose \eqref{existence_ult_tan_A_g1} and \eqref{existence_ult_tan} both hold. Then, it may be confirmed that \eqref{viable_case_z1} and \eqref{viable_case} \emph{both} do not hold, implying that there are no candidate barrier curves intersecting $G_-$ associated with the admissible set, giving a contradiction. Similarly, if \eqref{existence_ult_tan_M_g1} and \eqref{existence_ult_tan} both hold, then \eqref{MRPI_case_z1} and \eqref{unconquerable_case_MRPI} both do not hold. This completes the proof.
\end{pf}

\begin{rem}
	Lemma~\ref{lemma_3} says that barriers (if they exist) intersect either $g_1$ or $g_3$, but that there cannot exist two barriers, with one intersecting $g_1$ and the other $g_3$.
\end{rem}

\begin{lem}\label{lemma_4}
	Suppose \eqref{viable_case} (resp. \eqref{unconquerable_case_MRPI}) does not hold. Then, there cannot exist a barrier curve intersecting a point of ultimate tangentiality on $g_1$. Similarly, suppose \eqref{viable_case_z1} (resp. \eqref{MRPI_case_z1}) does not hold. Then, there cannot exist a barrier curve intersecting a point of ultimate tangentiality on $g_3$.
\end{lem}
\begin{pf}
	Suppose \eqref{viable_case} does not hold and there exists a barrier curve, associated with $\A$, intersecting a point of ultimate tangentiality on $g_1$. Then, \eqref{existence_ult_tan_A_g1} holds. From Lemma~\ref{lemma_3} we know that \eqref{existence_ult_tan} does not hold, thus:
	\[
		\bar{x}_2 \geq \frac{A_h\bar{x}_1}{A_h\bar{x}_1 + \gamma}.
	\]
	If we consider this inequality along with:
	\[
		A_m(A_h + \gamma)\bar{x}_2 + \gamma u_{\max} \leq A_mA_h,
	\]
	we see that:
	\[
		A_h(A_m + u_{\max})\bar{x}_1 + \gamma u_{\max} \leq A_mA_h,
	\]
	and thus \eqref{viable_case_z1} is violated, leading to a contradiction. The proof of the second statement follows along the same lines.
\end{pf}

\begin{rem}
	Lemma~\ref{lemma_4} is related to ``desperate'' cases, as we will shortly see. Intuitively, it says that if $\bar{x}_2$ is so small that the system becomes desperate, then it is \emph{not} possible to make $\bar{x}_1$ small enough for the system to not be desperate any more.
\end{rem}

\subsection{Characterisation of the sets}

\cite{deLara2016viableMOSQUITOS} characterised the Ross-Macdonald system as being either comfortable, desperate or viable. We show how the conditions for these cases can be obtained from the analysis in Subsections~\ref{subsec_ult_tan}-\ref{subsec_input_realisation}, and introduce a fourth case: the \emph{comfortable-viable case}.
\begin{definition}\label{def_cases_RM}
	We will say that the system \eqref{eq_RossMacdonald} is:
	\begin{itemize}
		\item \emph{comfortable}, provided $\mathcal{M} = \mathcal{A} = [0,\bar{x}_1]\times [0,\bar{x}_2]$,
		\item \emph{comfortable-viable}, provided $\mathcal{A} \subset [0,\bar{x}_1]\times [0,\bar{x}_2]$, $\mathcal{M} \subset [0,\bar{x}_1]\times [0,\bar{x}_2]$, and \eqref{eq_RossMacdonald} is neither desperate nor comfortable
		\item \emph{viable}, provided $\mathcal{A} \subset [0,\bar{x}_1]\times [0,\bar{x}_2]$, $\mathcal{M} = \{0\}\times \{0\}$, and \eqref{eq_RossMacdonald} is not desperate,
		\item \emph{desperate}, provided $\mathcal{M} = \mathcal{A} = \{0\}\times \{0\}$.
	\end{itemize} 
\end{definition}
The comfortable and desperate cases occur when the sets are trivial. The viable case occurs when $\M$ is trivial, but $\A$ is not. The new comfortable-viable case occurs when both sets are nontrivial.

\begin{rem}
	The set $\M$ may also be considered a ``comfortable'' set, in the sense that the state will always remain in it. We also note that, for the system \eqref{eq_RossMacdonald}, $\M$ is only nontrivial if $\A$ is, and so the comfortable-viable case is a special viable case. Hence our chosen name.
\end{rem} 

\begin{prop}\label{prop_cases}
	The system \eqref{eq_RossMacdonald} is:
	\begin{itemize}
		\item \emph{comfortable} $\Leftrightarrow$ both \eqref{existence_ult_tan_A_g1} and \eqref{existence_ult_tan} do not hold.
		\item \emph{comfortable-viable} $\Leftrightarrow$ either \eqref{existence_ult_tan} and \eqref{unconquerable_case_MRPI} hold; or \eqref{existence_ult_tan_M_g1} and \eqref{MRPI_case_z1} hold.
		\item \emph{viable} $\Leftrightarrow$ either \eqref{existence_ult_tan} and \eqref{viable_case} hold; or \eqref{existence_ult_tan_A_g1} and \eqref{viable_case_z1} hold.
		\item \emph{desperate} $\Leftrightarrow$ \eqref{existence_ult_tan} holds and \eqref{viable_case} does not $\Leftrightarrow$ \eqref{existence_ult_tan_A_g1} holds and \eqref{viable_case_z1} does not.
	\end{itemize}
\end{prop}
\begin{pf}
	\begin{itemize}
		\item The comfortable case occurs when the set $[0,\bar{x}_1]\times[0,\bar{x}_2]$ is invariant. Thus, we require 
		$$\max_{u\in U} L_fg_1(x,u)\leq 0\quad \mathrm{for~all~} x\in\{\bar{x}_1\}\times[0,\bar{x}_2]$$ 
		and 
		$$\max_{u\in U} L_fg_3(x,u)\leq 0\quad \mathrm{for~all~} x\in[0,\bar{x}_1]\times\{\bar{x}_2\}.$$	
		The result then follows from the fact that 
		$$\frac{A_mx_2}{A_mx_2 + u_{\min}}\leq \frac{A_m\bar{x}_2}{A_m\bar{x}_2 + u_{\min}}\quad \mathrm{for~all~} x_2\in[0,\bar{x}_2],$$ 
		and 
		$$\frac{A_hx_1}{A_hx_1 +\gamma}\leq \frac{A_h\bar{x}_1}{A_h\bar{x}_1 + \gamma}\quad \mathrm{for~all~} x_1\in[0,\bar{x}_1].$$
		\item The comfortable-viable case occurs when there exists an invariance barrier curve contained in $G_-$ that ends at a tangent point located either on $[0,\bar{x}_1[\times\{\bar{x}_2\}$ or on $\{\bar{x}_1\}\times[0,\bar{x}_2[$. (Recall, from Lemma~\ref{lemma_3}, that both cannot happen.) Thus, either there exists a tangent point on $g_1$ (condition \eqref{existence_ult_tan_M_g1}) along with a curve contained in $G_-$ (condition \eqref{MRPI_case_z1}); or the point is located on $g_3$ (condition \eqref{existence_ult_tan}) and the curve evolves backwards into $G_-$ (condition \eqref{unconquerable_case_MRPI}).
		\item Similar to the comfortable-viable case, the viable case occurs when there exit barrier curves, associated with the admissible set, contained in $G_-$ and ending at tangent points on $[0,\bar{x}_1[\times\{\bar{x}_2\}$ or $\{\bar{x}_1\}\times[0,\bar{x}_2[$.
		\item  The desperate case occurs when the system is not comfortable (thus there exists a tangent point on $[0,\bar{x}_1[\times\{\bar{x}_2\}$ or $\{\bar{x}_1\}\times[0,\bar{x}_2[$), and candidate barrier curves, associated with the admissible set, are \emph{not} contained in $G_-$. 
		For $g_1$, this is true if and only if \eqref{existence_ult_tan_A_g1} holds and \eqref{viable_case_z1} does not hold. From Lemma~\ref{lemma_4} we know that this is equivalent to \eqref{existence_ult_tan} holding and \eqref{viable_case} not holding.
	\end{itemize}
This completes the proof.
\end{pf}

\begin{figure}[h]
	\begin{center}
		\includegraphics[width=\linewidth,height=\linewidth,keepaspectratio]{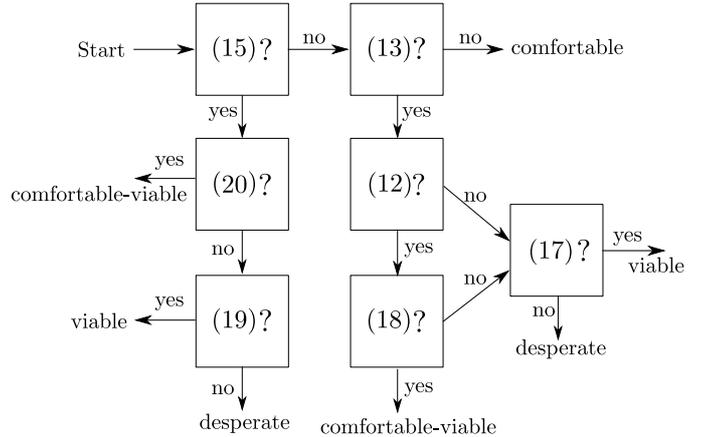} 
		\caption{A flow-diagram summarising the use of Proposition~\ref{prop_cases} to classify the system.} 
		\label{fig_flow_diagram}                             
	\end{center}                        
\end{figure}

Proposition~\ref{prop_cases} may be summarised with the flow-diagram shown in Figure~\ref{fig_flow_diagram}.

\begin{rem}
	In our characterisation of the sets we derive the same conditions for the comfortable, desperate and viable cases as those of \cite{deLara2016viableMOSQUITOS} where $\bar{x}_1 = 1$, but using different arguments that use the conditions of Theorem~\ref{thm1}. To now construct the sets (for the viable and comfortable-viable cases) one merely needs to integrate the dynamics and adjoint equations backwards from $z_1$ or $z_3$ using $\bar{u}$ until the resulting curve intersects a constraint.
\end{rem}

\section{Numerical Example}\label{sec_examples}

\begin{figure}[h]
	\begin{center}\
		\includegraphics[width=\linewidth,keepaspectratio]{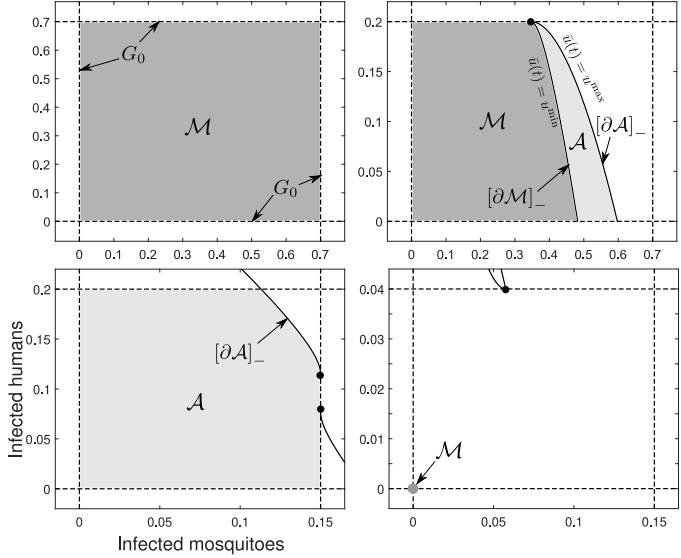} 
		\caption{The MRPI and admissible set for system \eqref{eq_RossMacdonald} under a comfortable ($\bar{x}_1 = \bar{x}_2 = 0.7$), comfortable-viable ($\bar{x}_1 = 0.7, \bar{x}_2 = 0.2$), viable ($\bar{x}_1 = 0.15, \bar{x}_2 = 0.2$) and desperate case ($\bar{x}_1 = 0.15, \bar{x}_2 = 0.04$).} 
		\label{fig:epidem_model}                             
	\end{center}                              
\end{figure}

We use the same model parameters as those estimated in Appendix C of \cite{deLara2016viableMOSQUITOS}, for a dengue outbreak in Cali, Colombia. The parameters are $A_m = 0.076608$, $A_h = 0.0722633, \gamma = 0.1$, with $u~\in~[u_{\min},u_{\max}]=[0.0333,0.05]$. We find the sets for varying bounds on the state, see Figure~\ref{fig:epidem_model}. We see that when the infection cap on both humans and mosquitoes are high ($\bar{x}_1 = 0.7$, $\bar{x}_2 = 0.7$) we are in a comfortable case, and the entire constraint set is robustly invariant. As we decrease the human infection cap, keeping the mosquito cap the same, we enter a comfortable-viable case ($\bar{x}_1 = 0.7$, $\bar{x}_2 = 0.2$), where there exists a nontrivial admissible set and MRPI. If we now decrease the mosquito cap we enter a viable case ($\bar{x}_1 = 0.15$, $\bar{x}_2 = 0.2$) with the barrier curve now intersecting $g_1$ instead of $g_3$, as it did previously. Further decreasing the bound on the infected humans ($\bar{x}_1 = 0.15$, $\bar{x}_2 = 0.04$) we enter a desperate case, and no sets exist. In the bottom two plots we also indicate the candidate barrier curves that evolve outside of $G_-$ from tangent points (black dots), and that must be ignored.

\section{Management of epidemics with the obtained sets}\label{sec_discussion}

We now describe how the results of our analysis may be used to maintain infection caps during epidemics. If the state is located in $\M$, then it is guaranteed that the infection caps will always be maintained regardless of fumigation strategy, as long as the minimal fumigation is always maintained. Thus, the set $\M$ is ``comfortable'', in the sense that only minimal resources are needed to maintain the infection caps. If the state is located in the interior of $\A$, then any fumigation strategy may be used for some period of time. However, if the state reaches the barrier of the admissible set, $\DAM$, then the only way to maintain the infection caps is to constantly use the maximal fumigation rate, $u_{\max}$. If the system is in the viable case, and the state is located at the infection caps but still in $\A$ (i.e. the state is at $\DAO$) then there is some freedom in the fumigation that may be used, but it is advisable to use $u_{\max}$. If the state is located outside both $\M$ and $\A$, then it is impossible to maintain the infection caps, and one is guaranteed to be violated in the future. The only way to prevent this is to relax the infection caps, or increase the maximal allowed fumigation.

To summarise, given the system parameters ($A_m$ and $A_h$), along with fumigation bounds ($u_{\max}$ and $u_{\min}$) and infection caps ($\bar{x}_1$ and $\bar{x}_2$), fumigation strategies should be implemented as follows. If the system is:
\begin{itemize}
	\item comfortable, and $x(t)\in[0,\bar{x}_1]\times [0,\bar{x}_2]$, let $u(t)=u_{\min}$.
	\item comfortable-viable and $x(t)\in\M$, let $u(t) = u_{\min}$.
	\item comfortable-viable and $x(t)\in \Int(\A)$, let $u(t) = u_{\min}$.
	\item comfortable-viable and $x(t)\in \DAM$, let $u(t) = u_{\max}$.
	\item viable and $x(t)\in \Int(\A)$, let $u(t) = u_{\min}$.
	\item viable and $x(t)\in \DAO \cup \DAM$, let $u(t) = u_{\max}$.
	\item any of the cases and $x(t)\in\Acomp$, then relax infection caps, or increase maximal fumigation.
	\item desperate, then relax the infection caps, or increase maximal allowed fumigation.
\end{itemize}

It is interesting to note that the more freedom there is in the fumigation rate, i.e. the larger the difference between $u_{\min}$ and $u_{\max}$, the more the sets differ, the ratio $\mathrm{Area}(\M)/\mathrm{Area}(\A)$ measuring the potential efficiency of the fumigation rate.

\section{Conclusion}\label{sec_conclusion}

We applied the theory of barriers to analyse the constrained Ross-Macdonald model of a vector-borne disease, aiming to maintain the proportion of infected humans and/or mosquitoes below specified caps. Building on the work by \cite{deLara2016viableMOSQUITOS}, we introduced the maximal robust positively invariant set to the study of epidemics, arguing that minimal resources are required to maintain infection caps when the state is in this set. We derived inequalities of the system parameters that may be used to classify the system into one of four classes, and we have constructed the sets for an example with various infection bounds. Future research could focus on analysing other epidemic systems such as the well-known ``compartmental'' models, that include interventions such as vaccination, \cite{HETHCOTE1973365}.

\bibliography{bib_epidemiology}             
\end{document}